\documentclass{amsart}
\usepackage{amssymb,euscript,amsmath, mathrsfs}
\usepackage{bm}
\usepackage{rotating}
\usepackage[table]{xcolor}

\makeindex

\newcounter{ENUM}
\newcommand{\itm}{\item}
\newenvironment{ilist}[1][0]{\renewcommand{\theENUM}{\roman{ENUM}}\renewcommand{\itm}{\addtocounter{ENUM}{1}\item[(\theENUM)]}\begin{itemize}\setcounter{ENUM}{#1}}{\end{itemize}}

\newenvironment{alist}[1][0]{\renewcommand{\theENUM}{\alph{ENUM}}\renewcommand{\itm}{\addtocounter{ENUM}{1}\item[(\theENUM)]}\begin{itemize}\setcounter{ENUM}{#1}}{\end{itemize}}

\newcommand{\margh}[1]{}

\def\risom{\overset{\sim}{\rightarrow}}

\input xy
\xyoption{all}
\CompileMatrices

\def\PP{{\mathbb P}}

\def\cS{{\mathcal S}}
\def\cL{{\mathcal L}}

\def\sO{{\mathscr O}}

\def\vp{\varphi}

\def\Pic{\operatorname{Pic}}

\def\exdm{\operatorname{expdim}}

\newtheorem{thm}{Theorem}[section]
\newtheorem{prop}[thm]{Proposition}
\newtheorem{lem}[thm]{Lemma}

\newtheorem{conj}[thm]{Conjecture}

\theoremstyle{definition}
\newtheorem{defn}[thm]{Definition}
\newtheorem{ques}[thm]{Question}

\newtheorem{ex}[thm]{Example}

\theoremstyle{remark}

\newtheorem{rem}[thm]{Remark}

\numberwithin{equation}{section}

\begin{document}
\title{An observation on $(-1)$-curves on rational surfaces}
\author{Olivia Dumitrescu}
\author{Brian Osserman}
\begin{abstract} We give an effective iterative characterization of the 
classes of (smooth, rational) $(-1)$-curves on the blowup of the projective 
plane at general points.
Such classes are characterized as having self-intersection
$-1$, arithmetic genus $0$, and intersecting every $(-1)$-curve of smaller
degree nonnegatively.
\end{abstract}

\thanks{Brian Osserman is partially supported
by a grant from the Simons Foundation \#279151.}

\maketitle

\section{Introduction}

Let $S$ be a smooth projective surface. We will call an irreducible curve
$C$ on $S$ a \textbf{$(-1)$-curve} if it is smooth and rational, and
$C \cdot C =-1$. The $(-1)$-curves are of fundamental importance in
various aspects of surface theory, but we are especially interested in
their role in the interpolation problem. Specifically, suppose that
we have $n$ points $P_1,\dots,P_n$ on $\PP^2$. Given also positive
integers $d,m_1,\dots,m_n$, we can ask the following basic question:

\begin{ques}\label{ques:interpolation} What is the dimension of the
space of homogeneous polynomials of degree $d$ in three variables
which vanish to order at least $m_i$ at $P_i$ for $i=1,\dots,n$?
\end{ques}

The naive expected dimension is given by 
\begin{equation}\label{eq:exp-dim} 
\exdm (d,m_1,\dots,m_n):= 
\max\left\{0,\binom{d+2}{2}-\sum_{i=1}^n \binom{m_i+1}{2}\right\}.
\end{equation}

The actual dimension is always at least the expected one, but
one quickly sees that the actual dimension may be larger, even when
the $P_i$ are general: the first example is the case $d=n=m_1=m_2=2$,
where the expected dimension is $0$, but in fact the space of polynomials
is $1$-dimensional, containing (the polynomial whose zero set is) the 
doubled line through $P_1$ and $P_2$. 

Now, let $S$ be the blowup of $\PP^2$ at $P_1,\dots,P_n$, and let
$H,E_1,\dots,E_n$ denote the classes of the hyperplane and the exceptional
divisors, respectively. Then our space of homogeneous polynomials can
be reinterpreted as $\Gamma(S,\sO(dH-m_1 E_1 - \dots - m_n E_n))$, and
we observe that the strict transform of a line through $P_1$ and $P_2$
is a $(-1)$-curve on $S$, which occurs twice in the base locus of
$\sO(2H-2E_1-2E_2)$. The fundamental conjecture in the field is that this 
phenomenon is the only one causing the expected dimension to differ from
the actual dimension:

\begin{conj} (Gimigliano-Harbourne-Hirschowitz \cite{gi3} \cite{ha4} 
\cite{hi2}) 
Assume the situation 
is as above, and the $P_i$ are general. Let 
$\cL=\sO(dH-m_1 E_1 - \dots - m_n E_n)$. Then we have that
$$\dim \Gamma(S,\cL)=\exdm (d,m_1,\dots,m_n)$$
unless there is some $(-1)$-curve $C$ in $S$ such that
$C \cdot \cL \leq -2$. 
\end{conj}

Thus, it becomes important to understand the $(-1)$-curves in surfaces
$S$ obtained from $\PP^2$ by blowing up general points. An effective
but iterative description of such curves as iterated images under Cremona 
transformations of exceptional divisors 
is exposited in Dolgachev (Corollary 1 of \cite{do1}).
The purpose of this note is to give an 
alternative description, which is also effective but iterative, but which
may be more satisfying conceptually. In essence, it says that a class
which numerically looks like the class of a $(-1)$-curve is in fact the
class of a $(-1)$-curve unless it contains a smaller-degree $(-1)$-curve
in its base locus. More specifically, our result is the following:

\begin{thm}\label{thm:main} Let $S$ be the blowup of $\PP^2$ at very
general points $P_1,\dots,P_n$. A divisor class on $S$ is the class of a
$(-1)$-curve if and only if either it is one of the $E_i$, or it is
of the form $dH-m_1 E_1 - \dots - m_n E_n$, with $d>0$, and $m_i \geq 0$
for all $i$, and the following conditions are satisfied:
\begin{ilist}
\itm (self-intersection $-1$) $d^2-\sum_i m_i^2=-1$;
\itm (arithmetic genus $0$) 
$\frac{(d-1)(d-2)}{2}-\sum_i \frac{m_i(m_i-1)}{2}=0$;
\itm for all $0< d' < d$, and all $(-1)$-curves $C'$ of degree $d'$ on $S$,
we have $C \cdot C' \geq 0$.
\end{ilist}
\end{thm}

Thus, we get a purely numerical (albeit inductive) criterion for
describing all the classes of $(-1)$-curves on $S$. If we consider a 
single class at a time, it suffices to assume that the $P_i$ are general;
we impose that they are very general in order to treat all classes 
simultaneously.

The three conditions of the theorem are clearly necessary. Conversely,
if conditions (i) and (ii) are satisfied, it is easy to see (Proposition
\ref{prop:num-equivs} below) that the class is effective, so the heart of
the matter is to show that if a curve in the class is not irreducible, it must
have negative product with a $(-1)$-curve having smaller degree.

\begin{rem} The last condition of the theorem is necessary: the lowest-degree 
example is in the case $d=5$, $n=10$, with the class 
$5H-3E_1-3E_2-E_3-\dots-E_{10}$. This satisfies conditions (i) and (ii),
but it is not the class of a $(-1)$-curve; indeed, it can be written as
the sum of the classes $H-E_1-E_2$ and $4H-2E_1-2E_2-E_3-\dots-E_{10}$,
which represent disjoint genus-$0$ and genus-$1$ curves, respectively.

In addition, the generality of points is necessary: for instance, if
$P_1,P_2,P_3$ are collinear, then the class $H-E_1-E_2$ satisfies the
hypotheses of the theorem, but is no longer the class of a $(-1)$-curve: 
rather, it is the class of the line through $P_1,P_2,P_3$ together with 
$E_3$.
\end{rem}

The proof of Theorem \ref{thm:main} is brief, and relies on the same
techniques (namely, a lemma of Noether allowing inductive application
of Cremona transformations) as the aforementioned criterion in \cite{do1}.

\subsection*{Acknowledgements} We would like to thank Joe Harris, who
suggested to us that Theorem \ref{thm:main} should hold, and Ciro
Ciliberto, who pointed out an omitted hypothesis in an earlier version
of Lemma \ref{lem:noether}.

\section{Proof}

Our proofs rely on two basic background elements: a numerical lemma of
Noether, and the theory of standard Cremona transformations of the plane.
Both of these are described in Dolgachev \cite{do1};
however, because their proofs are short and our hypotheses 
for Noether's lemma are slightly more general, 
we have elected to give a self-contained presentation,
nonetheless following the arguments in \cite{do1}.

Our first two results are purely numerical, and do not involve any
generality conditions on the points $P_i$.

\begin{prop}\label{prop:num-equivs} Let $C=dH-m_1 E_1 - \dots - m_n E_n$
be an arbitrary curve class on $S$. Then any two of the following 
conditions imply the remaining two conditions.

\begin{alist}
\itm ($C \cdot C = -1$) $d^2-\sum_i m_i^2=-1$;
\itm ($C$ has arithmetic genus $0$) 
$\frac{(d-1)(d-2)}{2}-\sum_i \frac{m_i(m_i-1)}{2}=0$;
\itm ($C \cdot (-K) = 1$) $3d-\sum_i m_i =1$;
\itm ($\exdm = 1$) $\frac{(d+2)(d+1)}{2}-\sum_i \frac{(m_i+1)m_i}{2}=1$.
\end{alist}
\end{prop}

\begin{proof} If $K$ is the canonical class of $S$, we have
$K=-3H+E_1 + \dots + E_n$, 
and by adjunction the arithmetic genus of $C$ is $\frac{2+C\cdot (C+K)}{2}$.
We therefore see that we can re-express all our conditions in terms of 
intersection
theory: (a) is $C \cdot C = -1$, (b) is $C \cdot C + C \cdot K = -2$,
(c) is $C \cdot K = -1$, and (d) is $C \cdot C= C \cdot K$. 
The proposition follows trivially.
\end{proof}

We now prove the aforementioned lemma of Noether; our proof follows 
\cite{do1}.
\margh{cite also [go]?}
We will only use the case $C \cdot C = -1$, but we have included the more
general statement because the proof is the same.

\begin{lem}\label{lem:noether} Let $C=dH-m_1 E_1 -\dots- m_n E_n$ be a
curve class with $d>0$ and all $m_i \geq 0$, and satisfying condition (b) 
of Proposition 
\ref{prop:num-equivs}, and also having $-2 \leq C \cdot C \leq 1$. If $d=1$,
suppose further that $C \cdot C < 0$. Then there exist $i_1<i_2<i_3$ such 
that
$$m_{i_1}+m_{i_2}+m_{i_3}>d.$$
\end{lem}

\begin{proof} Without loss of generality, we may suppose that we have
ordered our points so that $m_1 \geq m_2 \geq \dots \geq m_n$, and we
thus want simply to show that $m_1+m_2+m_3 > d$. If $d=1$, 
then the hypothesis $C \cdot C < 0$ means that we cannot have $m_1=m_2=0$
or $m_1=1$ and $m_2=0$, so the lemma holds in this case.
We therefore assume that $d>1$.

Now, for $1 \leq j \neq n$, define
$$q_j = \frac{\sum_{i=j}^n m_i^2}{\sum_{i=j}^n m_i}.$$
Then because $m_j \geq m_i$ for $i \geq j$, we have $m_j \geq q_j$ for any
$j$. Set $r_j=m_j-q_j$ for each $j$. For brevity, set $c=C \cdot C$,
and observe (using condition (b) of Proposition \ref{prop:num-equivs}) that 
we have the following expressions:
$$q_1=\frac{d^2-c}{3d-2-c}, \text{ and } q_i=q_{i-1}-r_{i-1} 
\frac{m_{i-1}}{m_i+\dots+m_n} \text{ for } i>1.$$
Condition (b) of Proposition \ref{prop:num-equivs} also gives us first that
$m_i<d$ for all $i$ (because $d>1$),
and second that we have $c=C \cdot C = -2-C \cdot K$, 
so that $3d-m_1 - \dots - m_n=c+2$. Because we have $m_i<d$ for $i=1,2$, 
it follows that we have
$$m_2 + \dots + m_n = 3d-m_1-c-2 \geq 2d-c-1,\text{ and }
m_3 + \dots + m_n = 3d-m_1-m_2-c-2 \geq d-c.$$
Again using $m_i < d$ for $i=1,2$, we thus find that 
$$q_2=q_1-r_1 \frac{m_1}{m_2+\dots+m_n} \geq q_1-r_1 \frac{d-1}{2d-c-1},$$
and
$$q_3=q_2-r_2 \frac{m_2}{m_3+\dots+m_n} \geq q_2-r_2 \frac{d-1}{d-c}.$$
We finally conclude that
\begin{align*} m_1+m_2+m_3 & \geq q_1+r_1+q_2+r_2+q_3 \geq
q_1+r_1+2 q_2 +r_2(1-\frac{d-1}{d-c}) \\
& \geq 3q_1 +r_1(1- 2 \frac{d-1}{2d-c-1})+r_2(1-\frac{d-1}{d-c}) \\
& \geq 3q_1 = 3\frac{d^2-c}{3d-2-c}>d,
\end{align*}
where we have used our hypothesis that $c \leq 1$ in the last two 
inequalities, and also that $d>1$ and $c \geq -2$ for the final inequality.
\end{proof}

For our purposes, we will be interested only in the Cremona transformations
which are standard Cremona transformations based in some three of our
points $P_i$. One can use the geometric Cremona transformation to realize
$S$ also as a blowup of $\PP^2$ in a different set of points; from this
point of view, the transformation induces a change of basis of $\Pic(S)$.
However, for our applications we will be interested instead in a more 
artificial construction, in which we keep our fixed basis of $\Pic(S)$, 
and use a numerical form of the Cremona transformation to construct an 
automorphism on $\Pic(S)$. This leads to the following definition.

\begin{defn}\label{def:cremona} Given $1 \leq i_1<i_2<i_3 \leq n$, the
$(i_1,i_2,i_3)$-Cremona transformation on $S$ is the automorphism of
$\Pic(S)$ which sends a class $dH-m_1 E_1 - \dots - m_n E_n$ to 
$d'H-m'_1 E_1 - \dots - m'_n E_n$, where $d'=2d-m_{i_1}-m_{i_2}-m_{i_3}$,
$m'_i=m_i$ for $i\neq i_1,i_2,i_3$, and
$m_{i_1}'=d-m_{i_2}-m_{i_3}$, 
$m_{i_2}'=d-m_{i_1}-m_{i_3}$, and
$m_{i_3}'=d-m_{i_1}-m_{i_2}$.
\end{defn}

Note that this map is an automorphism because it is visibly a homomorphism,
and it is an involution. Geometrically speaking, when no three of the $P_i$
are collinear, this map is motivated
by replacing the exceptional divisors $E_{i_1},E_{i_2},E_{i_3}$ with
the $(-1)$-curves $H-E_{i_2}-E_{i_3}$, $H-E_{i_1}-E_{i_3}$, and 
$H-E_{i_1}-E_{i_2}$, and changing the $P_j$ for $j \neq i_1,i_2,i_3$
by the induced geometric Cremona transformation. Although our purely
numerical definition is a distinct construction, in the below proposition
we leverage the connection between our definition and its
geometric motivation to show that when the $P_i$ are general, our
definition has good behavior.

The following proposition is essentially one half of Theorem 1 of
\cite{do1}.

\begin{prop}\label{prop:cremona-props} Assume that $P_1,\dots,P_r$
are very general in $\PP^2$. Then Cremona transformations on $S$ preserve
the following:
\begin{ilist}
\itm Intersection products;
\itm The canonical class;
\itm Effective classes.
\end{ilist}
\end{prop}

If we consider any given class, it suffices to assume that the $P_i$ are
general; we assume they are very general in order to treat all classes
simultaneously.

\begin{proof} That Cremona transformations preserve intersection products
and the canonical class is an elementary calculation which holds regardless
of generality of the $P_i$.
For the preservation of effective classes, we will more precisely show
the following: if we fix a class $C=dH-m_1 E_1 - \dots - m_n E_n$, and
$i_1<i_2<i_3$, for general choice of the $P_i$, we will have that $C$ is
effective if and only if the $(i_1,i_2,i_3)$-Cremona image of $C$ is
effective. In fact, since the $(i_1,i_2,i_3)$-Cremona map is an
involution on curve classes, it is enough to show that the image of an
effective class is effective. Fix any choice of $P_{i_1},P_{i_2},P_{i_3}$
which are not collinear, and let $U \subseteq \PP^2$ be the complement of
the three lines through pairs of the $P_{i_1},P_{i_2},P_{i_3}$, so that
$U^{n-3}$ parametrizes
the choices of the remaining $P_j$ such that none of the $P_j$ are
collinear with any two of $P_{i_1},P_{i_2},P_{i_3}$. Thus, the geometric
Cremona transformation based at $P_{i_1},P_{i_2},P_{i_3}$ can be thought
of as an automorphism $\vp:U \risom U$, inducing an automorphism 
$\vp^j:U^j \risom U^j$. If $U' \subseteq U^j$ is the complement of the
pairwise diagonals, then $\vp^j$ maps $U'$ into itself, inducing an
automorphism $\vp':U' \to U'$. Let $\cS$ be the surface over $U'$ obtained 
from $\PP^2\times U'$ by blowing up $P_{i_1},P_{i_2},P_{i_3}$ (considered as
constant sections) together with the $j$ disjoint universal sections over
$U'$. Then via the geometry of the Cremona transformation, the fiber
$S$ of $\cS$ over a given choice of $(P_j)_j \in U'$ can simultaneously be 
realized as the fiber over $(\vp(P_j))_j$, and the resulting change of
basis of $\Pic(S)$ is precisely the (numerical) 
$(i_1,i_2,i_3)$-Cremona transformation.

Now, a class 
$C=dH-m_1 E_1 -\dots-m_n E_n$ induces a line bundle on $\cS$, and since
$\cS$ is flat and projective over $U'$, 
by the semicontinuity theorem either $C$ is effective over all points of $U'$,
or there is a nonempty open subset $U'' \subseteq U'$ over which $C$ is not
effective. But we see by construction that $C$ is effective over a given
point $(P_j)_j \in U'$ if and only if its $(i_1,i_2,i_3)$-Cremona 
transformation is effective over its Cremona image $\vp'((P_j)_j)$. 
Thus, if $C$ is
effective over a general point of $U'$, it is effective over all of $U'$,
and its $(i_1,i_2,i_3)$-Cremona transformation is likewise effective over 
all of $U'$.
On the other hand, if $C$ fails to be effective over $U'' \subseteq U'$,
then provided we choose $(P_j)_j \in U''$, the statement of the proposition
will still be satisfied.
\end{proof}

\begin{ex} Note that we need some generality hypothesis on the $P_i$
in order for Proposition \ref{prop:cremona-props} to hold, as otherwise we 
could have for instance 
$P_1,\dots,P_6$ all lying on a conic, so that $2H-E_1-\dots-E_6$ is
effective, but the Cremona transformation at $(1,2,3)$ gives
$H-E_4-E_5-E_6$, which is not effective unless $P_4,P_5,P_6$ are
collinear.
\end{ex}

We are now ready to prove our main theorem.

\begin{proof}[Proof of Theorem \ref{thm:main}]
We first observe that a curve class is the class
of a $(-1)$-curve if and only if it satisfies conditions (i) and (ii) of
the theorem, and cannot be written as a sum of effective classes. Certainly,
these conditions are necessary,
and conversely, if conditions (i) and (ii) are satisfied, Proposition 
\ref{prop:num-equivs} gives us that the expected dimension is positive, so 
the class is effective. If it is not a sum of effective classes, any
representative is an integral curve, hence a $(-1)$-curve. Consequently,
Proposition \ref{prop:cremona-props} also implies that Cremona 
transformations preserve classes of $(-1)$-curves.

Now, we see that the class of a $(-1)$-curve must satisfy conditions 
(i)-(iii) of the theorem, and we will prove the converse by induction on 
$d$. For the base case is $d=0$, we prove a stronger statement: even 
without the hypothesis $m_i \geq 0$, we note that conditions (i) and (ii) 
already imply
that $C$ is the class of one of the $E_i$, which is indeed a
$(-1)$-curve. Now, suppose that $d>0$, and we know the theorem for all 
$d'<d$. Given a class $C=dH-m_1 E_1 - \dots - m_n E_n$ satisfying conditions
(i) and (ii) of the statement, by our above observations, it suffices to
show that if $C$ is a sum of effective classes, then there is a $(-1)$-curve 
of smaller degree intersecting $C$ negatively. In fact, we observe that the 
smaller-degree condition is automatic: since the $(-1)$-curve would 
necessarily occur in the base locus of $C$, its degree is at most $d$.
If its degree is exactly $d$, the difference is a positive combination
of the $E_i$,
but we see immediately from Proposition \ref{prop:num-equivs} that we
cannot have two curve classes satisfying conditions (i) and (ii) which
differ by a positive combination of the $E_i$. 

Now, according to Lemma \ref{lem:noether}, there exist $i_1<i_2<i_3$
such that $m_{i_1}+m_{i_2}+m_{i_3}>d$. We therefore apply an 
$(i_1,i_2,i_3)$-Cremona transformation, yielding a 
curve class $C'=d' H - m_1' E_1 - \dots - m_n' E_n$, 
where $d'=2d-m_{i_1}-m_{i_2}-m_{i_3}<d$. 
Then by 
Proposition \ref{prop:cremona-props}, we have that $C'$ still satisfies
conditions (i) and (ii) of the theorem, and is still effective. In
particular, $d' \geq 0$. If $C=C_1+C_2$ is a sum of effective classes,
applying the same Cremona transformation to $C_1$ and $C_2$, we write
$C'=C'_1+C'_2$, with $C'_1$ and $C'_2$ still effective by Proposition
\ref{prop:cremona-props}, so in particular $C'$ is not the class of a 
$(-1)$-curve. First note that if $d'=0$, then by the base case we have that 
$C'= \pm E_i$, contradicting that $C'$ is a not a $(-1)$-curve.
On the other hand, if $d'>0$, we conclude that 
there is a $(-1)$-curve $C''$ of degree $d''<d'$ such that $C'' \cdot C'<0$.
If $m'_i \geq 0$ for all $i$, this follows from the induction hypothesis, 
while if some $m'_i<0$, we set $C''=E_i$.
Applying the same Cremona transformation to $C''$ produces a $(-1)$-curve
(of distinct class from $C$) intersecting $C$ negatively, as desired.
\end{proof}

Thus, if we argue iteratively instead of inductively, we see that if we
start with a non-$(-1)$-curve class, the process always ends with one
of the $m_i$ going negative.


\bibliographystyle{amsalpha}
\bibliography{gen}

\end{document}